\documentstyle[12pt]{article}

\begin{document}

\centerline{\bf NONCLASSICAL REPRESENTATIONS OF THE NONSTANDARD}
\centerline{\bf DEFORMATIONS $U'_q({\rm so}_n)$, $U_q({\rm iso}_n)$ AND
$U'_q({\rm so}_{n,1})$}

\bigskip
\bigskip
\bigskip
\centerline{N.~Z.~Iorgov \  and  \ A.~U.~Klimyk}
\centerline{Institute for Theoretical Physics,
Kiev 252143, Ukraine}

\bigskip
\bigskip

\begin{abstract}
The aim of this paper is to announce the results on irreducible
nonclassical type representations of the nonstandard $q$-deformations
$U'_q({\rm so}_n)$, $U_q({\rm iso}_n)$ and $U'_q({\rm so}_{n,1})$
of the universal enveloping algebras of the Lie algebras
${\rm so}(n,{\bf C})$, ${\rm iso}_n$ and ${\rm so}_{n,1}$ when $q$ is
a real number (the algebra $U'_q({\rm so}_{n,1})$ is a real form of the
algebra $U'_q({\rm so}_{n+1})$). These representations are characterized by
the properties that they are singular at the point $q=1$.

\end{abstract}

\bigskip
\bigskip

\centerline{\bf 1. Introduction}
\bigskip

    Quantum orthogonal groups, quantum Lorentz group and their
corresponding quantum algebras are of special interest for modern
physics [1-3]. M. Jimbo [4] and V. Drinfeld [5] defined $q$-deformations
(quantum algebras) $U_q(g)$ for all simple complex Lie algebras $g$
by means of Cartan subalgebras and root subspaces (see also [6]).
Reshetikhin,
Takhtajan and Faddeev [7] defined quantum algebras $U_q(g)$ in terms
of the universal $R$-matrix. However, these approaches do not
give a satisfactory presentation of the quantum algebra
$U_q({\rm so}(n,{\bf C}))$ from a viewpoint of some problems in quantum
physics and representation theory.
When considering representations of the quantum
groups $SO_q(n+1)$ and $SO_q(n,1)$
we are interested in reducing them onto the quantum subgroup $SO_q(n)$.
This reduction would give the analogue of the Gel'fand-Tsetlin basis for
these representations. However, definitions of quantum algebras
mentioned above do not allow the inclusions
$U_q({\rm so}(n+1, {\bf C}))\supset U_q({\rm so}(n,{\bf C}))$
and $U_q({\rm so}_{n,1})\supset U_q({\rm so}_n)$. To be able to exploit such
reductions we have to consider $q$-deformations of the Lie algebra
${\rm so}(n+1,{\bf C})$ defined in terms of the generators
$I_{k,k-1}=E_{k,k-1}-E_{k-1,k}$  (where $E_{is}$ is the matrix with
elements $(E_{is})_{rt}=\delta _{ir} \delta _{st})$ rather than
by means of Cartan subalgebras and root elements.
To construct such deformations we have to deform trilinear relations
for elements $I_{k,k-1}$ instead of Serre's relations (as in the case of
Jimbo's quantum algebras). As a result, we obtain the associative algebra
which will be denoted as $U'_q({\rm so}(n,{\bf C})).$

       These $q$-deformations were first constructed  in [8].
They permit one to construct the reductions of $U'_q({\rm so}_{n,1})$ and
$U'_q({\rm so}_{n+1})$ onto $U'_q({\rm so}_n)$.
The $q$-deformed algebra $U'_q({\rm so}(n,{\bf C}))$ leads for $n=3$ to
the $q$-deformed algebra $U'_q({\rm so}(3,{\bf C}))$ defined by
D. Fairlie [9]. The cyclically symmetric algebra, similar to Fairlie's
one, was also considered somewhat earlier by Odesskii [10].
The algebra $U'_q({\rm so}(3,{\bf C}))$ allows us to construct the noncompact
quantum algebra $U'_q({\rm so}_{2,1})$. The algebra
$U'_q({\rm so}(4,{\bf C}))$ is a $q$-deformation of the Lie
algebra ${\rm so}(4,{\bf C})$ given by means of usual bilinear
commutation relations between the elements $I_{ji}$, $1\le i<j\le 4$.
In the case of the classical
Lie algebra ${\rm so}(4,{\bf C})$ one has
${\rm so}(4,{\bf C})={\rm so}(3,{\bf C})+{\rm so}(3,{\bf C})$,
while in the case of our $q$-deformation $U'_q({\rm so}(4,{\bf C}))$
this is not the case.

In the classical case, the imbedding $SO(n)\subset SU(n)$
(and its infinitesimal analogue) is of great importance for nuclear
physics and in the theory of Riemannian spaces. It is well known
that in the framework of Drinfeld--Jimbo quantum groups and algebras
one cannot construct the corresponding imbedding. The algebra
$U'_q({\rm so}(n,{\bf C}))$ allows to define such an imbedding [11],
that is, it is possible to define the imbedding
$U'_q({\rm so}(n,{\bf C}))\subset U_q({\rm sl}_n)$,
where $U_q({\rm sl}_n)$ is the Drinfeld-Jimbo quantum algebra.

As a disadvantage of the algebra $U'_q({\rm so}(n,{\bf C}))$ we have
to mention the difficulties with Hopf algebra structure. Nevertheless,
$U'_q({\rm so}(n,{\bf C}))$ turns out to be a coideal in
$U_q({\rm sl}_n).$

Finite dimensional irreducible representations of the algebra
$U'_q({\rm so}(n,{\bf C}))$ were constructed in [8]. The formulas
of action of the generators of $U'_q({\rm so}(n,{\bf C}))$ upon the
basis (which is a $q$-analogue of the Gel'fand--Tsetlin basis) are
given there. A proof of these formulas and some their corrections were
given in [12]. However,
finite dimensional irreducible representations described in [8] and [12]
are representations of the classical type. They are $q$-deformations of the
corresponding irreducible representations of the Lie algebra
${\rm so}(n,{\bf C})$, that is, at $q\to 1$ they turn into representations
of ${\rm so}(n,{\bf C})$.

The algebra $U'_q({\rm so}(n,{\bf C}))$ has other classes of finite
dimensional irreducible representations which have no classical analogue.
These representations are singular at the limit $q\to 1$. One of the aims
of this paper is to describe these representations of
$U'_q({\rm so}(n,{\bf C}))$. Note that the description of these
representations for the algebra $U'_q({\rm so}(3,{\bf C}))$ is given in
[13]. A classification of irreducible $*$-representations of real forms
of the algebra $U'_q({\rm so}(3,{\bf C}))$ is given in [14].

There exists an algebra, closely related to the algebra
$U'_q({\rm so}(n,{\bf C}))$, which is a $q$-deformation of the universal
enveloping algebra $U({\rm iso}_n)$ of the Lie algebra ${\rm iso}_n$ of the
Euclidean group $ISO(n)$ (see [15]). It is denoted as
$U_q({\rm iso}_n)$. Irreducible representations of the classical type of
the algebra $U_q({\rm iso}_n)$ were described in [15]. A proof of the
corresponding formulas was given in [16]. However, the algebra
$U_q({\rm iso}_n)$, $q\in {\bf R}$, has irreducible representations of the
nonclassical type. A description of these representations is the second
aim of this paper. Note that the description of these representations for
$U_q({\rm iso}_2)$ is given in [17]. A classification of irreducible
$*$-representations of $U_q({\rm iso}_2)$ is obtained in [18].
The last aim of this paper is to describe irreducible representations of
nonclassical type of the algebra $U'_q({\rm so}_{n,1})$ which is a real
form of the algebra $U'_q({\rm so}(n+1,{\bf C}))$. Representations of the
classical type of this algebra are described in [8] and [20].

We assume throughout the paper that $q$ is a fixed positive number.
Thus, we give formulas for representations for these values of $q$. However,
these representations can be considered for any values of $q$ not
coinciding with a root of unity. For this we have to treat appropriately
square roots in formulas for representations or to rescale basis
vector in such a way that formulas for representations would not
contain square roots.

For convenience, we denote the Lie algebra ${\rm so}(n,{\bf C})$ by
${\rm so}_n$ and the algebra $U'_q({\rm so}(n,{\bf C}))$ by
$U'_q({\rm so}_n)$.
\bigskip

\centerline{\bf 2. The $q$-deformed algebras $U'_q({\rm so}_n)$ and
$U_q({\rm iso}_n)$}
\medskip

In our approach [8] to the
$q$-deformation of the algebras $U({\rm so}_n)$ we define the $q$-deformed
algebras $U'_q({\rm so}(n,{\bf C}))$ as the associate algebra (with a unit)
generated by the elements $I_{i,i-1}$, $i=2,3,...,n$
satisfying the defining relations
$$
I_{i,i-1}I^2_{i-1,i-2}-(q+q^{-1})I_{i-1,i-2}I_{i,i-1}I_{i-1,i-2}
+I^2_{i-1,i-2}I_{i,i-1}=-I_{i,i-1},     \eqno (1)
$$
$$
I^2_{i,i-1}I_{i-1,i-2}-(q+q^{-1})I_{i,i-1}I_{i-1,i-2}
I_{i,i-1}+I_{i-1,i-2}I^2_{i,i-1}=-I_{i-1,i-2},   \eqno(2)
$$
$$
[I_{i,i-1},I_{j,j-1}]=0,\ \ \ \vert i-j\vert >1,    \eqno (3)
$$
where [.,.] denotes usual commutator. Obviously, in the limit $q\to 1$
formulas (1)--(3) give the relations
$$
I_{i,i-1}I^2_{i-1,i-2}-2I_{i-1,i-2}I_{i,i-1}I_{i-1,i-2}
+I^2_{i-1,i-2}I_{i,i-1}=-I_{i,i-1},
$$
$$
I^2_{i,i-1}I_{i-1,i-2}-2I_{i,i-1}I_{i-1,i-2}
I_{i,i-1}+I_{i-1,i-2}I^2_{i,i-1}=-I_{i-1,i-2},
$$
$$
[I_{i,i-1},I_{j,j-1}]=0,\ \ \ \vert i-j\vert >1,
$$
defining the universal enveloping
algebra $U({\rm so}_n)$. Note also that relations
(1) and (2) principally differ from the $q$-deformed Serre relations
in the approach of Jimbo [4] and Drinfeld [5] to quantum orthogonal
algebras by a presence of nonzero right hand side and by possibility
of the reduction
$$
U'_q({\rm so}_n)\supset U'_q({\rm so}_{n-1})\supset
\cdots \supset U'_q({\rm so}_3).
$$
Recall that in the standard Jimbo--Drinfeld approach to the
definition of quantum algebras, the algebras
$U_q({\rm so}_{2m})$ and the algebras
$U_q({\rm so}_{2m+1})$ are distinct series of quantum algebras
which are constructed independently of each other.

   Various real forms of the algebras $U'_q({\rm so}_n)$ are
obtained by imposing corresponding $*$-structures (antilinear
antiautomorphisms). The compact real form $U'_q({\rm so}(n))$ is defined
by the $*$-structure
$$I^*_{i,i-1}=-I_{i,i-1},\ \ \  i=2,3,...,n.  $$
The noncompact $q$-deformed algebras $U'_q({\rm so}_{p,r})$ where $r=n-p$
are singled out respectively by means of the $*$-structures
$$I^*_{i,i-1}=-I_{i,i-1},\ \  i\ne p+1,\ \  i\le n, \ \ \
I^*_{p+1,p}=I_{p+1,p}.  $$
Among the noncompact real $q$-algebras $U'_q({\rm so}_{p,r}),$
the algebras $U'_q({\rm so}_{n-1,1})$
(a $q$-analogue of the Lorentz algebras) are of special interest.
\medskip

We also define the algebra $U_q({\rm iso}_n)$ which is a
nonstandard deformation of the universal enveloping algebra of the
Lie algebra ${\rm iso}_n$ of the Euclidean Lie group $ISO(n)$.
It is the associative algebra (with a unit) generated by the elements
$I_{21}$, $I_{32}, \cdots ,I_{n,n-1}, T_n$ such that the elements
$I_{21}$, $I_{32}, \cdots ,I_{n,n-1}$ satisfy the defining relations
of the subalgebra $U'_q({\rm so}_n)$ and the additional
defining relations are
$$
I^2_{n,n-1}T_n -(q+q^{-1}) I_{n,n-1}T_nI_{n,n-1}
+T_nI^2_{n,n-1}=-T_n, $$
$$
T_n^2I_{n,n-1} -(q+q^{-1}) T_nI_{n,n-1}T_n
+I_{n,n-1}T^2_n=0, $$
$$
[I_{k,k-1},T_n]\equiv I_{k,k-1}T_n-T_nI_{k,k-1}=0\ \ \ {\rm if}\ \ \
k<n
$$
(see [15]).
If $q=1$, then these relations define the classical algebra
$U({\rm iso}_n)$.

Let us note that the defining relations for $U_q({\rm iso}_n)$ can be
expressed by bilinear relations [16]. As an example, we consider the algebra
$U_q({\rm iso}_2)$ (see [16] and [18]).
This algebra is generated by the elements $I_{21}$ and
$T_2$.
Setting $T_1=[I_{21},T_2]_q\equiv q^{1/2}I_{21}T_2
-q^{-1/2}T_2I_{21}$, we obtain from the two defining relations the
bilinear relations
$$
[I_{21},T_2]_q=T_1,\ \ \ \
[T_1,I_{21}]_q=T_2,\ \ \ \
[T_2,T_1]_q=0
$$
which are also defining relations for $U_q({\rm iso}_2)$.
Note that the elements $T_1$ and $T_2$ corresponding to infinitesimal
generators of shifts in the Lie algebra ${\rm iso}_2$ do not commute
(they $q$-commute,  that is, $q^{1/2}T_2T_1-q^{-1/2}T_1T_2=0$).
A similar picture we have for the algebra $U_q({\rm iso}_n)$.
\medskip

\centerline{\bf 3. Finite dimensional classical type
representations of $U'_q({\rm so}_n)$}
\medskip

In this section we describe
(in the framework of a $q$-analogue of Gel'fand--Tsetlin formalism)
irreducible finite dimensional
representations of the algebras $U'_q({\rm so}_{n})$, $n \ge 3$,
which are $q$-deformations of the finite dimensional irreducible
representations of the Lie algebra ${\rm so}_n$.
They are given by sets ${\bf m}_{n}$
consisting of $\{ {n/2}\}$ numbers $m_{1,n}, m_{2,n},..., m_{\left
\{ {n/2}\right \} ,n}$ (here $\{ {n/2}\}$ denotes integral part
of ${n/2}$) which are all integral or all half-integral and
satisfy the dominance conditions
$$
m_{1,2p+1}\ge m_{2,2p+1}\ge ... \ge
m_{p,2p+1}\ge 0 ,                                  \eqno(4)
$$
$$
m_{1,2p}\ge m_{2,2p}\ge ... \ge m_{p-1,2p}\ge |m_{p,2p}|
                                                    \eqno(5)
$$
for $n=2p+1$ and $n=2p$, respectively.
These representations are denoted by $T_{{\bf m}_n}$.
For a basis in a representation space we take the $q$-analogue of
Gel'fand--Tsetlin basis which is obtained by successive reduction of
the representation $T_{{\bf m}_n}$ to the subalgebras
$U'_q({\rm so}_{n-1})$, $U'_q({\rm so}_{n-2})$, $\cdots$, $U'_q({\rm so}_3)$,
$U'_q({\rm so}_2):=U({\rm so}_2)$.
As in the classical case, its elements are labelled by Gel'fand--Tsetlin
tableaux
$$
  \{\xi_{n} \}
\equiv \left\{ \matrix{ {\bf m}_{n} \cr {\bf m}_{n-1} \cr \dots
\cr {\bf m}_{2}  }
 \right\}
\equiv \{ {\bf m}_{n},\xi_{n-1}\}\equiv \{{\bf m}_{n} ,
{\bf m}_{n-1} ,\xi_{n-2}\} ,
                                                       \eqno(6)
$$
where the components of ${\bf m}_{n}$ and ${\bf m}_{n-1}$ satisfy the
"betweenness" conditions
$$
m_{1,2p+1}\ge m_{1,2p}\ge m_{2,2p+1} \ge m_{2,2p} \ge ...
\ge m_{p,2p+1} \ge m_{p,2p} \ge -m_{p,2p+1}  ,
                                                     \eqno(7)
$$
$$
m_{1,2p}\ge m_{1,2p-1}\ge m_{2,2p} \ge m_{2,2p-1} \ge ...
\ge m_{p-1,2p-1} \ge \vert m_{p,2p} \vert .
                                                      \eqno(8)
$$
The basis element defined by tableau $\{\xi_{n} \}$ is denoted
as $\vert \{\xi_{n} \} \rangle $ or simply as $\vert
\xi_{n} \rangle $.

It is convenient to introduce the so-called $l$-coordinates
$$
l_{j,2p+1}=m_{j,2p+1}+p-j+1,  \qquad
                      l_{j,2p}=m_{j,2p}+p-j ,        \eqno(9)
$$
for the numbers $m_{i,k}$.
In particular, $l_{1,3}=m_{1,3}+1$ and $l_{1,2}=m_{1,2}$.
The operator $T_{{\bf m}_n}(I_{2p+1,2p})$ of the representation
$T_{{\bf m}_n}$ of $U'_q({\rm so}_{n})$ acts upon Gel'fand--Tsetlin
basis elements, labelled by (6), by the formula
$$
T_{{\bf m}_n}(I_{2p+1,2p})
| \xi_n\rangle =
\sum^p_{j=1} \frac{ A^j_{2p}(\xi_n)}
{q^{l_{j,2p}}+q^{-l_{j,2p}} }
            \vert (\xi_n)^{+j}_{2p}\rangle -
\sum^p_{j=1} \frac{A^j_{2p}((\xi_n)^{-j}_{2p})}
{q^{l_{j,2p}}+q^{-l_{j,2p}}}
|(\xi_n)^{-j}_{2p}\rangle                         \eqno(10)
$$
and the operator $T_{{\bf m}_n}(I_{2p,2p-1})$ of the representation
$T_{{\bf m}_n}$ acts as
$$
T_{{\bf m}_n}(I_{2p,2p-1})\vert \xi_n\rangle=
\sum^{p-1}_{j=1} \frac{B^j_{2p-1}(\xi_n)}
{[2 l_{j,2p-1}-1][l_{j,2p-1}]}
\vert (\xi_n)^{+j}_{2p-1} \rangle
$$
$$
-\sum^{p-1}_{j=1}\frac {B^j_{2p-1}((\xi_n)^{-j}_{2p-1})}
{[2 l_{j,2p-1}-1][l_{j,2p-1}-1]}
\vert (\xi_n)^{-j}_{2p-1}\rangle
+ {\rm i}\, C_{2p-1}(\xi_n)
\vert \xi_n \rangle .                                 \eqno(11)
$$
In these formulas, $(\xi_n)^{\pm j}_{k}$ means the tableau (6)
in which $j$-th component $m_{j,k}$ in ${\bf m}_k$ is replaced
by $m_{j,k}\pm 1$. The coefficients
$A^j_{2p},  $ $B^j_{2p-1},$ $C_{2p-1}$ in (10) and (11) are given
by the expressions
$$
A^j_{2p}(\xi_n) =
\left( \frac{\prod_{i=1}^p [l_{i,2p+1}+l_{j,2p}] [l_{i,2p+1}-l_{j,2p}-1]
\prod_{i=1}^{p-1} [l_{i,2p-1}+l_{j,2p}] [l_{i,2p-1}-l_{j,2p}-1]}
{\prod_{i\ne j}^p [l_{i,2p}+l_{j,2p}][l_{i,2p}-l_{j,2p}]
[l_{i,2p}+l_{j,2p}+1][l_{i,2p}-l_{j,2p}-1]} \right)^{1/2},  \eqno (12)
$$
and
$$
B^j_{2p-1}(\xi_n)=\left( \frac{\prod_{i=1}^p
[l_{i,2p}+l_{j,2p-1}] [l_{i,2p}-l_{j,2p-1}] \prod_{i=1}^{p-1}
[l_{i,2p-2}+l_{j,2p-1}] [l_{i,2p-2}-l_{j,2p-1}]}
{\prod_{i\ne j}^{p-1}
[l_{i,2p-1}{+}l_{j,2p-1}][l_{i,2p-1}{-}l_{j,2p{-}1}]
[l_{i,2p-1}{+}l_{j,2p-1}{-}1][l_{i,2p-1}{-}l_{j,2p-1}{-}1]}  \right) ^{1/2} ,
\eqno(13)
$$
$$
C_{2p-1}(\xi_n) =\frac{ \prod_{s=1}^p [ l_{s,2p} ]
\prod_{s=1}^{p-1} [ l_{s,2p-2} ]}
{\prod_{s=1}^{p-1} [l_{s,2p-1}] [l_{s,2p-1} - 1] } ,   \eqno(14)
$$
where numbers in square brackets mean $q$-numbers defined by
$$
[a]:= \frac {q^a-q^{-a}}{q-q^{-1}}.
$$
In particular,
$$
T_{{\bf m}_n}(I_{3,2})\vert \xi_n\rangle=
\frac{1}{q^{m_{1,2}}+q^{-m_{1,2}}} ( ([l_{1,3}+m_{1,2}]
[l_{1,3}-m_{1,2}-1])^{1/2} \vert (\xi_n )^{+1}_2\rangle -   $$
$$
-([l_{1,3}+m_{1,2}-1][l_{1,3}-m_{1,2}])^{1/2}
\vert (\xi_n )^{-1}_2\rangle  ) ,  $$
$$
T_{{\bf m}_n}(I_{2,1})\vert \xi_n\rangle= {\rm i}[m_{1,2}]
\vert \xi_n\rangle ,
$$

It is seen from (9) that $C_{2p-1}$ in (14) identically
vanishes if $m_{p,2p}\equiv l_{p,2p}=0$.

A proof of the fact that formulas (10)-(14) indeed determine
a representation of $U'_q({\rm so}_n)$ is given in [12].
\medskip

\centerline{\bf 4. Finite dimensional nonclassical type
representations of $U'_q({\rm so}_n)$}
\medskip

The representations of the previous section are called representations
of the classical type, since under the limit $q\to 1$ the operators
$T_{{\bf m}_n}(I_{j,j-1})$ turn into the corresponding operators
$T_{{\bf m}_n}(I_{j,j-1})$ for irreducible finite dimensional
representations with highest weights ${\bf m}_n$ of the Lie algebra
${\rm so}_n$.

The algebra $U'_q({\rm so}_n)$ also has irreducible finite dimensional
representations $T$ of nonclassical type, that is, such that the operators
$T(I_{j,j-1})$ have no classical limit $q\to 1$.
They are given by sets $\epsilon := (\epsilon _2,\epsilon _3,\cdots ,
\epsilon _n)$, $\epsilon _i=\pm 1$, and by sets
${\bf m}_{n}$ consisting of $\{ {n/2}\}$ {\bf half-integral} numbers
$m_{1,n}, m_{2,n}, \cdots , m_{\{ n/2\} ,n}$
(here $\{ {n/2}\}$ denotes integral part of ${n/2}$) that satisfy the
dominance conditions
$$
m_{1,2p+1}\ge m_{2,2p+1}\ge ... \ge
m_{p,2p+1}\ge 1/2 ,                                  \eqno(15)
$$
$$
m_{1,2p}\ge m_{2,2p}\ge ... \ge m_{p-1,2p}\ge m_{p,2p}\ge 1/2
                                                    \eqno(16)
$$
for $n=2p+1$ and $n=2p$, respectively.
These representations are denoted by $T_{\epsilon,{\bf m}_n}$.

For a basis in the representation space we use the analogue of the
basis of the previous section. Its elements are
labelled by tableaux
$$
  \{\xi_{n} \}
\equiv \left\{ \matrix{ {\bf m}_{n} \cr {\bf m}_{n-1} \cr \dots
\cr {\bf m}_{2}  }
 \right\}
\equiv \{ {\bf m}_{n},\xi_{n-1}\}\equiv \{{\bf m}_{n} ,
{\bf m}_{n-1} ,\xi_{n-2}\} ,
                                                       \eqno(17)
$$
where the components of ${\bf m}_{n}$ and ${\bf m}_{n-1}$ satisfy the
"betweenness" conditions
$$
m_{1,2p+1}\ge m_{1,2p}\ge m_{2,2p+1} \ge m_{2,2p} \ge ...
\ge m_{p,2p+1} \ge m_{p,2p} \ge 1/2  ,
                                                     \eqno(18)
$$
$$
m_{1,2p}\ge m_{1,2p-1}\ge m_{2,2p} \ge m_{2,2p-1} \ge ...
\ge m_{p-1,2p-1} \ge m_{p,2p}  .
                                                      \eqno(19)
$$
The basis element defined by tableau $\{\xi_{n} \}$ is denoted
as $\vert \xi_{n} \rangle $.

As in the previous section, it is convenient to introduce the
$l$-coordinates
$$
l_{j,2p+1}=m_{j,2p+1}+p-j+1,  \qquad
                      l_{j,2p}=m_{j,2p}+p-j .        \eqno(20)
$$
The operator $T_{\epsilon,{\bf m}_n}(I_{2p+1,2p})$ of the representation
$T_{\epsilon,{\bf m}_n}$ of $U_q({\rm so}_{n})$ acts upon our
basis elements, labelled by (17), by the formulas
$$
T_{\epsilon,{\bf m}_n}(I_{2p+1,2p})
| \xi_n\rangle =
\delta_{m_{p,2p},1/2}\, \frac{\epsilon _{2p+1}}{q^{1/2}-q^{-1/2}} D_{2p}
(\xi _n) | \xi_n\rangle +
$$
$$
+\sum^p_{j=1} \frac{ A^j_{2p}(\xi_n)}
{q^{l_{j,2p}}-q^{-l_{j,2p}} }
            \vert (\xi_n)^{+j}_{2p}\rangle -
\sum^p_{j=1} \frac{A^j_{2p}((\xi_n)^{-j}_{2p})}
{q^{l_{j,2p}}-q^{-l_{j,2p}}}
|(\xi_n)^{-j}_{2p}\rangle    ,                     \eqno(21)
$$
where the summation in the last sum must be from 1 to $p-1$ if
$m_{p,2p}=1/2$,
and the operator $T_{{\bf m}_n}(I_{2p,2p-1})$ of the representation
$T_{{\bf m}_n}$ acts as
$$
T_{\epsilon ,{\bf m}_n}(I_{2p,2p-1})\vert \xi_n\rangle=
\sum^{p-1}_{j=1} \frac{B^j_{2p-1}(\xi_n)}
{[2 l_{j,2p-1}-1][l_{j,2p-1}]_+}
\vert (\xi_n)^{+j}_{2p-1} \rangle -
$$
$$
-\sum^{p-1}_{j=1}\frac {B^j_{2p-1}((\xi_n)^{-j}_{2p-1})}
{[2 l_{j,2p-1}-1][l_{j,2p-1}-1]_+}
\vert (\xi_n)^{-j}_{2p-1}\rangle
+ \epsilon _{2p} {\hat C}_{2p-1}(\xi_n)
\vert \xi_n \rangle ,                                 \eqno(22)
$$
where
$$
[a]_+=\frac {q^a+q^{-a}}{q-q^{-1}}.
$$
In these formulas, $(\xi_n)^{\pm j}_{k}$ means the tableau (17)
in which $j$-th component $m_{j,k}$ in ${\bf m}_k$ is replaced
by $m_{j,k}\pm 1.$ Matrix elements
$A^j_{2p}$ and $B^j_{2p-1}$ are given by the same formulas
as in (10) and (11) (that is, by the formulas (12) and (13)) and
$$
{\hat C}_{2p-1}(\xi_n) = {
\prod_{s=1}^p [ l_{s,2p} ]_+
\prod_{s=1}^{p-1} [ l_{s,2p-2} ]_+ \over
   \prod_{s=1}^{p-1} [l_{s,2p-1}]_+ [l_{s,2p-1} - 1]_+ } .   \eqno(23)
$$
$$
D_{2p} (\xi _n)=
\frac{\prod_{i=1}^p
[l_{i,2p+1}-\frac 12 ] \prod_{i=1}^{p-1} [l_{i,2p-1}-\frac 12 ] }
{\prod_{i=1}^{p-1}
[l_{i,2p}+\frac 12 ] [l_{i,2p}-\frac 12 ] } . \eqno (24)
$$

For the operators $T_{\epsilon,{\bf m}_n}(I_{3,2})$ and
$T_{\epsilon,{\bf m}_n}(I_{2,1})$ we have
$$
T_{\epsilon ,{\bf m}_n}(I_{3,2})\vert \xi_n\rangle=
\frac{1}{q^{m_{1,2}}-q^{-m_{1,2}}} ( ([l_{1,3}+m_{1,2}]
[l_{1,3}-m_{1,2}-1])^{1/2} \vert (\xi_n )^{+1}_2\rangle -  $$
$$
([l_{1,3}+m_{1,2}-1][l_{1,3}-m_{1,2}])^{1/2}
\vert (\xi_n )^{-1}_2\rangle )
$$
if $m_{1,2}\ne \frac 12$,
$$
T_{\epsilon ,{\bf m}_n}(I_{3,2})\vert \xi_n\rangle =
\frac{1}{q^{1/2} - q^{-1/2}} (\epsilon _3 [l_{1,3}-1/2]
\vert (\xi_n )\rangle +
([l_{1,3}+1/2][l_{1,3}-3/2])^{1/2}
\vert (\xi_n )^{+1}_2\rangle )
$$
if $m_{1,2} = \frac 12$, and
$$
T_{\epsilon ,{\bf m}_n}(I_{2,1})\vert \xi_n\rangle= \epsilon _2 [m_{1,2}]_+
\vert \xi_n\rangle .
$$

The fact that the above operators $T_{\epsilon ,{\bf m}_n}(I_{k,k-1})$
satisfy the defining relations (1)--(3) of the algebra $U'_q({\rm so}_n)$
is proved in the following way. We take the formulas (10)--(14) for
the classical type representations $T_{{\bf m}_n}$ of $U'_q({\rm so}_n)$
with half-integral $m_{i,n}$ and replace there every
$m_{j,2p+1}$ by $m_{j,2p+1}-{\rm i}\pi /2h$,
every $m_{j,2p}$, $j\ne p$, by $m_{j,2p}-{\rm i}\pi /2h$
and $m_{p,2p}$ by $m_{p,2p}-\epsilon _2 \epsilon _4\cdots \epsilon _{2p}
{\rm i}\pi /2h$, where each $\epsilon _{2s}$ is equal
to $+1$ or $-1$ and $h$ is defined by $q=e^h$.
Repeating almost word by word the reasoning of the paper [12],
we prove that the operators given by
formulas (10)--(14) satisfy the defining relations (1)--(3)
of the algebra $U'_q({\rm so}_n)$ after this replacement. Therefore,
these operators determine a representation of $U'_q({\rm so}_n)$. We
denote this representation by $T'_{{\bf m}_n}$.
After a simple rescaling, the operators
$T'_{{\bf m}_n}(I_{k,k-1})$ take the form
$$
T'_{{\bf m}_n}(I_{2p+1,2p})
| \xi_n\rangle =
\sum^p_{j=1} \frac{ A^j_{2p}(\xi_n)}
{q^{l_{j,2p}}-q^{-l_{j,2p}} }
            \vert (\xi_n)^{+j}_{2p}\rangle -
\sum^p_{j=1} \frac{A^j_{2p}((\xi_n)^{-j}_{2p})}
{q^{l_{j,2p}}-q^{-l_{j,2p}}}
|(\xi_n)^{-j}_{2p}\rangle    ,
$$
$$
T'_{{\bf m}_n}(I_{2p,2p-1})\vert \xi_n\rangle=
\sum^{p-1}_{j=1} \frac{B^j_{2p-1}(\xi_n)}
{[2 l_{j,2p-1}-1][l_{j,2p-1}]_+}
\vert (\xi_n)^{+j}_{2p-1} \rangle -
$$
$$
-\sum^{p-1}_{j=1}\frac {B^j_{2p-1}((\xi_n)^{-j}_{2p-1})}
{[2 l_{j,2p-1}-1][l_{j,2p-1}-1]_+}
\vert (\xi_n)^{-j}_{2p-1}\rangle
+ \epsilon _{2p} {\hat C}_{2p-1}(\xi_n)
\vert \xi_n \rangle ,
$$
where $A^j_{2p}$, $B^j_{2p-1}$ and ${\hat C}_{2p-1}$ are such as in
the formulas (21) and (22).
The representations
$T'_{{\bf m}_n}$ are reducible. We decompose these representations into
subrepresentations in the following way.
We fix $p$ ($p=1,2,\cdots ,\{ (n-1)/2\}$)
and decompose the space ${\cal H}$ of the
representation $T'_{{\bf m}_n}$ into direct sum of two
subspaces ${\cal H}_{\epsilon _{2p+1}}$,
$\epsilon _{2p+1}=\pm 1$, spanned by the basis vectors
$$
| \xi_n \rangle _{\epsilon _{2p+1}} =| \xi_n \rangle -\epsilon _{2p+1}
| \xi_n' \rangle , \ \ \ \ \ m_{p,2p}\ge 1/2,
$$
respectively,
where $| \xi_n' \rangle$ is obtained from $| \xi_n \rangle$ by replacement
of $m_{p,2p}$ by $-m_{p,2p}$.
A direct verification shows that two subspaces
${\cal H}_{\epsilon _{2p+1}}$ are
invariant with respect to all the operators $T'_{{\bf m}_n}(I_{k,k-1})$.
Now we take the subspaces ${\cal H}_{\epsilon _{2p+1}}$ and repeat the
same procedure for some $s$, $s\ne p$, and decompose each of these
subspaces into two invariant subspaces. Continuing this procedure further
we decompose the representation space ${\cal H}$ into a direct sum
of $2^{\{ (n-1)/2\} }$ invariant subspaces. The operators
$T'_{{\bf m}_n}(I_{k,k-1})$ act upon these subspaces by the formulas
(21) and (22). We denote the corresponding subrepresentations on these
subspaces by $T_{\epsilon ,{\bf m}_n}$. The above reasoning shows that
the operators
$T_{\epsilon ,{\bf m}_n}(I_{k,k-1})$ satisfy the defining relations
(1)--(3) of the algebra $U'_q({\rm so}_n)$.
\medskip

\noindent
{\bf Theorem 1.} {\it The representations $T_{\epsilon ,{\bf m}_n}$
are irreducible.
The representations $T_{\epsilon ,{\bf m}_n}$ and
$T_{\epsilon ',{\bf m}'_n}$ are pairwise nonequivalent for
$(\epsilon ,{\bf m}_n)\ne (\epsilon ',{\bf m}'_n)$.
For any admissable $(\epsilon ,{\bf m}_n)$ and ${\bf m}'_n$
the representations $T_{\epsilon ,{\bf m}_n}$ and $T_{{\bf m}'_n}$ are
pairwise nonequivalent.}
\medskip

The algebra $U'_q({\rm so}_n)$ has non-trivial one-dimensional
representations. They are special cases of the representations of the
nonclassical type. They are described as follows.

Let $\epsilon := (\epsilon _2,\epsilon _3,\cdots ,
\epsilon _n)$, $\epsilon _i=\pm 1$, and let
${\bf m}_{n}=(m_{1,n}, m_{2,n}, \cdots , m_{ \{ n/2 \} ,n})
=(\frac 12 , \frac 12 ,\cdots ,\frac 12 )$. Then the corresponding
representations $T_{\epsilon ,{\bf m}_n} $ are one-dimensional
and are given by the formulas
$$
T_{\epsilon ,{\bf m}_n}(I_{k+1,k})| \xi_n\rangle =
\frac{\epsilon _{k+1}}{q^{1/2}-q^{-1/2}}.
$$
Thus, to every $\epsilon := (\epsilon _2,\epsilon _3,\cdots ,
\epsilon _n)$, $\epsilon _i=\pm 1$, there corresponds a one-dimensional
representation of $U'_q({\rm so}_n)$.
\medskip

\centerline{\bf 5. Definition of representations of $U'_q({\rm so}_{n,1})$
and $U_q({\rm iso}_{n})$}
\medskip

Let us recall that we assume that $q$ is a positive number.
We give the following definition of infinite dimensional
representations of the algebras $U'_q({\rm so}_{n,1})$ and
$U_q({\rm iso}_n)$ (we denote these algebras by ${\cal A}$).
It is a homomorphism $R:
{\cal A}\to {\cal L}({\cal H})$
of ${\cal A}$ to the space ${\cal L}({\cal H})$ of linear
operators (bounded or unbounded) on a Hilbert space ${\cal H}$ such that
\medskip

(a) operators $R(a)$, $a\in {\cal A}$, are defined on an
invariant everywhere dense subspace ${\cal D}\subset {\cal H}$;
\smallskip

(b) $R\downarrow U'_q({\rm so}_n)$ decomposes into a direct sum of
irreducible finite dimensional representations of $U'_q({\rm so}_n)$
(with finite multiplicities if $R$ is irreducible);
\smallskip

(c) subspaces of irreducible representations of $U'_q({\rm so}_n)$
belong to ${\cal D}$.
\medskip

Two infinite dimensional irreducible representations $R$ and $R'$
of ${\cal A}$ on spaces ${\cal H}$ and ${\cal H}'$, respectively,
are called (algebraically) equivalent if there exists an everywhere
dense invariant subspaces $V\subset {\cal D}$ and $V'\subset {\cal D}'$
and a one-to-one linear operator $A: V\to V'$ such that
$AR(a)v=R'(a)Av$ for all $a\in {\cal A}$ and $v\in V$.

Remark that our definition of infinite dimensional representations of
$U'_q({\rm so}_{n,1})$ and $U_q({\rm iso}_n)$ corresponds to the definition
of Harish-Chandra modules for the pairs $({\rm so}_{n,1},{\rm so}_n)$
and $({\rm iso}_n,{\rm so}_n)$, respectively. Thus, modules
determined by representations of the above definition can be called
$q$-Harish-Chandra modules of the pairs
$(U'_q({\rm so}_{n,1}),U'_q({\rm so}_n))$
and $(U_q({\rm iso}_n),U'_q({\rm so}_n))$, respectively.
\medskip

\centerline{\bf 6. Classical type representations of $U_q({\rm iso}_{n})$}
\medskip

There are the following classes of irreducible representations of
$U_q({\rm iso}_n)$:
\medskip

(a) Finite dimensional irreducible representations $R$ of
$U'_q({\rm so}_n)$. They are irreducible representations of
$U_q({\rm iso}_n)$ with $R(T_n)=0$.
\smallskip

(b) Infinite dimensional irreducible representations
of the classical type.
\smallskip

(c) Infinite dimensional irreducible representations
of the nonclassical type.
\smallskip

Let us describe representations of class (b). They are given by
non-zero complex
parameter $\lambda$ and by numbers ${\bf m}=(m_{2,n+1},m_{3,n+2},\cdots ,
m_{\{ (n+1)/2\} ,n+1})$ describing irreducible representations of
the classical type of the subalgebra $U'_q({\rm so}_{n-1})$ (see
[15] and [16]).
We denote the corresponding representations of $U_q({\rm iso}_n)$ by
$R_{\lambda {\bf m}}$.

In order to describe the space of the representation
$R_{\lambda {\bf m}}$ we note that
$$
R_{\lambda {\bf m}}\downarrow U'_q({\rm so}_{n}) =
\bigoplus\nolimits _{{\bf m}_n}
T_{{\bf m}_n},\ \ \ {\bf m}_n=(m_{1,n},\cdots ,m_{\{ n/2 \} ,n}) ,
\eqno (26)
$$
where the summation is over all irreducible representations $T_{{\bf m}_n}$
of $U'_q({\rm so}_{n})$ which contain the irreducible representation
of $U'_q({\rm so}_{n-1})$ given by the numbers ${\bf m}$,
that is, such that
$$
m_{1,n}\ge m_{2,n+1}\ge m_{2,n}\ge m_{3,n+1}\ge m_{3,n}\ge \cdots .
$$
The carrier space ${\hat {\cal H}}_{\bf m}$ of the representation
$R_{\lambda {\bf m}}$ decomposes as
${\hat {\cal H}}_{\bf m}=\bigoplus _{{\bf m}_n} {\cal H}_{{\bf m}_n}$, where
the summation is such as in (26) and ${\cal H}_{{\bf m}_n}$ are the
subspaces, where the representations $T_{{\bf m}_n}$ of
$U'_q({\rm so}_{n})$ are realized. We choose the Gel'fand--Tsetlin basis
in every subspace ${\cal H}_{{\bf m}_n}$. The set of all these
Gel'fand--Tsetlin bases gives a basis of the space ${\hat {\cal H}}_{\bf m}$.
We denote the basis elements by $| {\bf m}_n,M\rangle$, where $M$ are
the corresponding Gel'fand--Tsetlin tableaux. The numbers $m_{ij}$ from
$| {\bf m}_n,M\rangle$ determine the numbers $l_{ij}$ as in section 3.
The numbers $m_{i,n+1}$ determine the numbers
$$
l_{i,2k+1}=m_{i,2k+1}+k-i+1,\ \ n=2k,\ \ \ \
l_{i,2k}=m_{i,2k}+k-i,\ \  n=2k-1.
$$
The operators
$R_{\lambda {\bf m}}(I_{i,i-1})$ are given by formulas of the classical type
representations of the algebra $U'_q({\rm so}_{n})$ given in section 3.
For the operators $R_{\lambda {\bf m}}(T_{2k})$ and
$R_{\lambda {\bf m}}(T_{2k-1})$ we have the expressions
$$
R_{\lambda {\bf m}}(T_{2k})|{\bf m}_{2k}, M\rangle =\lambda
\sum _{j=1} ^k \frac{{\tilde A}_{2k}^j ({\bf m}_{2k},M)}
{q^{l_{j,2k}}+q^{-l_{j,2k}}}
|{\bf m}^{+j}_{2k}, M\rangle + $$
$$\qquad\qquad\qquad
+\lambda
\sum _{j=1} ^k \frac{{\tilde A}_{2k}^j ({\bf m}^{-j}_{2k},M)}
{q^{l_{j,2k}}+q^{-l_{j,2k}}}
| m^{-j}_{2k},M\rangle , \eqno (27) $$
$$
R_{\lambda {\bf m}}(T_{2k-1})|{\bf m}_{2k-1}, M\rangle =\lambda
\sum _{j=1} ^{k-1} \frac {{\tilde B}_{2k-1}^j ({\bf m}_{2k-1},M)}
{[2l_{j,2k-1}-1][l_{j,2k-1}]}
| {\bf m}^{+j}_{2k-1},M\rangle + $$
$$
+\lambda
\sum _{j=1} ^{k-1} \frac{{\tilde B}_{2k-1}^j ({\bf m}^{-j}_{2k-1},M)}
{[2l_{j,2k-1}-1][l_{j,2k-1}-1]}
| {\bf m}^{-j}_{2k-1},M\rangle
+\lambda {\tilde C}_{2k-1}({\bf m}_{2k-1},M) |{\bf m}_{2k-1},M\rangle ,
\eqno (28)
$$
where ${\bf m}^{\pm j}_r$ means the set of numbers ${\bf m}_r$
with $m_{j,r}$
replaced by $m_{j,r}\pm 1$, respectively, the coefficients are given by
$$
{\tilde A}^j_{2k}({\bf m}_{2k},M)= \left(
\frac {\prod _{i=2}^{k} [l_{i,2k+1}+l_{j,2k}][l_{i,2k+1}-l_{j,2k}-1]}
{\prod _{i\ne j} [l_{i,2k}+l_{j,2k}][l_{i,2k}-l_{j,2k}]} \times \right.
$$
$$\qquad\qquad\qquad\qquad \left. \times
\frac {\prod _{i=1}^{k-1} [l_{i,2k-1}+l_{j,2k}][l_{i,2k-1}-l_{j,2k}-1]}
{\prod _{i\ne j} [l_{i,2k}+l_{j,2k}+1][l_{i,2k}-l_{j,2k}-1]}
\right) ^{1/2} ,
$$
$$
{\tilde B}^j_{2k-1}({\bf m}_{2k-1},M)= \left(
\frac {\prod _{i=2}^k [l_{i,2k}+l_{j,2k-1}][l_{i,2k}-l_{j,2k-1}]}
{\prod _{i\ne j} [l_{i,2k-1}+l_{j,2k-1}][l_{i,2k-1}-l_{j,2k-1}]}\times \right.
$$
$$\qquad\qquad\qquad  \left. \times
\frac {\prod _{i=1}^{k-1} [l_{i,2k-2}+l_{j,2k-1}][l_{i,2k-2}-l_{j,2k-1}]}
{\prod _{i\ne j} [l_{i,2k-1}+l_{j,2k-1}-1]
[l_{i,2k-1}-l_{j,2k-1}-1]} \right) ^{1/2} ,
$$
$$
{\tilde C}_{2k-1}(M)=
\frac {\prod _{s=2}^k [l_{s,2k}]\prod _{s=1}^{k-1} [l_{s,2k-2}]}
{\prod _{s=1}^{k-1} [l_{s,2k-1}][l_{s,2k-1}-1]} .
$$
These formulas were given in [16] (see also [15]).

The representations $ R_{\lambda {\bf m}}$, $\lambda \ne 0$, are irreducible.
The representations $R_{\lambda {\bf m}}$ and
$R_{\mu {\bf m}'}$ are equivalent if and only if ${\bf m}={\bf m}'$ and
$\lambda =\pm \mu$. The operator $R_{\lambda {\bf m}}(T_n)$ is bounded.
\medskip

\centerline{\bf 7. Nonclassical type representations of $U_q({\rm iso}_{n})$}
\medskip

Now we describe irreducible representations of nonclassical type (that is,
representations $R$ for which there exists no limit $q\to 1$
for the operators $R(T_n)$
and $R(I_{i,i-1})$). These representations are given by
$\epsilon :=(\epsilon _2,\epsilon _3,\cdots ,\epsilon _{n+1})$,
non-zero complex
parameter $\lambda$ and by numbers ${\bf m}=(m_{2,n+1},m_{3,n+2},\cdots ,
m_{\{ (n+1)/2\} ,n+1})$, $m_{2,n+1}\ge m_{3,n+2} \ge \cdots \ge
m_{\{ (n+1)/2\} ,n+1}\ge 1/2$,
describing irreducible representations of
the nonclassical type of the subalgebra $U'_q({\rm so}_{n-1})$ (see
section 4). We denote the corresponding representations of
$U_q({\rm iso}_n)$ by $R_{\epsilon ,\lambda ,{\bf m}}$.

In order to describe the space of the representation
$R_{\epsilon,\lambda ,{\bf m}}$ we note that
$$
R_{\epsilon,\lambda ,{\bf m}}\downarrow U'_q({\rm so}_{n}) =
\bigoplus _{{\bf m}_n}
T_{\epsilon ',{\bf m}_n},\ \ \ {\bf m}_n=(m_{1,n},\cdots ,m_{\{ n/2 \} ,n}) ,
\eqno (29)
$$
where $\epsilon '=(\epsilon _2,\cdots ,\epsilon _{n})$
is the part of the set $\epsilon$,
the summation is over all irreducible nonclassical type
representations $T_{\epsilon ',{\bf m}_n}$
of $U'_q({\rm so}_{n})$ for which
the components of ${\bf m}_{n}$ satisfy the
"betweenness" conditions
$$
m_{1,2k+1}\ge m_{1,2k}\ge m_{2,2k+1} \ge m_{2,2k} \ge ...
\ge m_{k,2k+1} \ge m_{k,2k} \ge 1/2 \ \ {\rm if}\ \ n=2k ,
$$
$$
m_{1,2k}\ge m_{1,2k-1}\ge m_{2,2k} \ge m_{2,2k-1} \ge ...
\ge m_{k-1,2k-1} \ge m_{k,2k} \ \ {\rm if}\ \ n=2k-1.
$$
The carrier space ${\hat {\cal H}}_{\epsilon,\bf m}$ of the representation
$R_{\epsilon,\lambda ,{\bf m}}$ decomposes as
${\hat {\cal H}}_{\epsilon, \bf m}=\bigoplus _{{\bf m}_n}
{\cal H}_{\epsilon ',{\bf m}_n}$, where
the summation is such as in (29) and ${\cal H}_{\epsilon',{\bf m}_n}$ are the
subspaces, where the representations $T_{\epsilon',{\bf m}_n}$ of
$U'_q({\rm so}_{n})$ are realized. We choose a basis in every subspace
${\cal H}_{\epsilon ',{\bf m}_n}$ as in section 4. The set of all these
bases gives a basis of the space ${\hat {\cal H}}_{\epsilon,\bf m}$.
We denote the basis elements by $| {\bf m}_n,M\rangle$, where $M$ are
the corresponding tableaux. The numbers $m_{ij}$ from
$| {\bf m}_n,M\rangle$ determine the numbers $l_{ij}$ as in section 3.
The numbers $m_{i,n+1}$ determine the numbers
$l_{i,n+1}$ as in section 6.
The operators $R_{\epsilon,\lambda ,{\bf m}}(I_{i,i-1})$ are given by
formulas of the nonclassical type
representations of the algebra $U'_q({\rm so}_{n})$ from
section 4. For the operators $R_{\epsilon,\lambda ,{\bf m}}(T_{2k})$ and
$R_{\epsilon,\lambda ,{\bf m}}(T_{2k-1})$ we have the expressions
$$
R_{\epsilon,\lambda ,{\bf m}}(T_{2k-1})|{\bf m}_{2k-1}, M\rangle =\lambda
\sum _{j=1} ^{k-1} \frac{ {\tilde B}_{2k-1}^j ({\bf m}_{2k-1},M)}
{[2l_{j,2k-1}-1][l_{j,2k-1}]_+}
| {\bf m}^{+j}_{2k-1},M\rangle + $$
$$
+\lambda
\sum _{j=1} ^{k-1} \frac{ {\tilde B}_{2k-1}^j ({\bf m}^{-j}_{2k-1},M)}
{[2l_{j,2k-1}-1][l_{j,2k-1}-1]_+}
| {\bf m}^{-j}_{2k-1},M\rangle
+{\rm i}\epsilon _{2k}
\lambda {\hat C}_{2k-1}({\bf m}_{2k-1},M) |{\bf m}_{2k-1},M\rangle ,
$$
$$
R_{\epsilon,\lambda ,{\bf m}}(T_{2k})|{\bf m}_{2k}, M\rangle ={\rm i}\lambda
\delta_{m_{p,2p},1/2}\, \frac{\epsilon _{2k+1}}{q^{1/2}-q^{-1/2}} D_{2k}
|{\bf m}_{2k}, M\rangle +     $$
$$
+\lambda \sum _{j=1} ^k \frac{{\tilde A}_{2k}^j ({\bf m}_{2k},M)}
{q^{l_{j,2k}}-q^{-l_{j,2k}}}
|{\bf m}^{+j}_{2k}, M\rangle +
\lambda
\sum _{j=1} ^k \frac{{\tilde A}_{2k}^j ({\bf m}^{-j}_{2k},M)}
{q^{l_{j,2k}}-q^{-l_{j,2k}}}
| {\bf m}^{-j}_{2k},M\rangle , $$
where the summation in the last sum must be from 1 to $k-1$ if
$m_{k,2k}=1/2$, and
${\tilde A}_{2k}^j $, ${\tilde B}_{2k-1}^j $ are such as in
(27) and (28), and
$$
{\hat C}_{2k-1}(M)=
\frac {\prod _{s=2}^k [l_{s,2k}]_+\prod _{s=1}^{k-1} [l_{s,2k-2}]_+}
{\prod _{s=1}^{k-1} [l_{s,2k-1}]_+[l_{s,2k-1}-1]_+} ,
$$
$$
D_{2k}=
\frac{\prod_{i=2}^k
[l_{i,2k+1}-\frac 12 ] \prod_{i=1}^{k-1} [l_{i,2k-1}-\frac 12 ] }
{\prod_{i=1}^{k-1}
[l_{i,2k}+\frac 12 ] [l_{i,2k}-\frac 12 ] } . \eqno (30)
$$

\noindent
{\bf Theorem 2.} {\it The representations $R_{\epsilon,\lambda ,{\bf m}}$
are irreducible. The representations $R_{\epsilon,\lambda ,{\bf m}}$ and
$R_{\epsilon ',\lambda ',{\bf m}'}$ are equivalent if and only if
$\epsilon =\epsilon '$, ${\bf m}={\bf m}'$ and $\lambda =\pm \lambda '$.
The operators $R_{\epsilon ,\lambda ,{\bf m}}(T_n)$ are bounded.
The representation $R_{\epsilon,\lambda ,{\bf m}}$ is equivalent to no
of the representations $R_{\lambda ',{\bf m}'}$ of section 6.}
\medskip

\centerline{\bf 8. Classical type representations of
$U'_q({\rm so}_{n,1})$}
\bigskip

Let us describe the principal series representations of the algebra
$U'_q({\rm so}_{n,1})$. They are given by a complex
parameter $c$ and by numbers ${\bf m}=(m_{2,n+1},m_{3,n+1},\cdots ,
m_{\{ (n+1)/2\} ,n+1})$ describing irreducible representations of
the classical type of the subalgebra $U'_q({\rm so}_{n-1})$ (see [8]
and [20]). We denote the corresponding representations of
$U_q({\rm so}_{n,1})$ by $R_{c,{\bf m}}$.

In order to describe the space of the representation
$R_{c,{\bf m}}$ we note that
$$
R_{c,{\bf m}}\downarrow U'_q({\rm so}_{n}) =
\bigoplus\nolimits _{{\bf m}_n}
T_{{\bf m}_n},\ \ \ {\bf m}_n=(m_{1,n},\cdots ,m_{\{ n/2 \},n}) ,
\eqno (31)
$$
where the summation is over all irreducible representations $T_{{\bf m}_n}$
of $U'_q({\rm so}_{n})$ which contain the irreducible representation
of $U'_q({\rm so}_{n-1})$ given by the numbers ${\bf m}$, that is,
$$
m_{1,n}\ge m_{2,n+1}\ge m_{2,n}\ge m_{3,n+1}\ge \cdots .
$$
Thus, the carrier space ${\hat {\cal H}}_{\bf m}$ of the representation
$R_{c,{\bf m}}$ decomposes as
${\hat {\cal H}}_{\bf m}=\bigoplus _{{\bf m}_n} {\cal H}_{{\bf m}_n}$, where
the summation is such as in (31) and ${\cal H}_{{\bf m}_n}$ are the
subspaces, where the representations $T_{{\bf m}_n}$ of
$U'_q({\rm so}_{n})$ are realized. We choose the Gel'fand--Tsetlin basis
in every subspace ${\cal H}_{{\bf m}_n}$. The set of all these
Gel'fand--Tsetlin bases gives a basis of the space ${\hat {\cal H}}_{\bf m}$.
We denote the basis elements by $| {\bf m}_n,M\rangle$, where $M$ are
the corresponding Gel'fand--Tsetlin tableaux. The numbers $m_{ij}$ from
$| {\bf m}_n,M\rangle$ determine the numbers $l_{ij}$ as in section 3.
The numbers $m_{i,n+1}$ determine the numbers
$$
l_{i,2k+1}=m_{i,2k+1}+k-i+1,\ \ n=2k,\ \ \ \
l_{i,2k}=m_{i,2k}+k-i,\ \  n=2k-1.
$$
The operators
$R_{c,{\bf m}}(I_{i+1,i})$, $i<n$, are given by formulas
of representations of the algebra $U'_q({\rm so}_{n})$ given in section 3.
For the action of the operators $R_{c,{\bf m}}(I_{2k+1,2k})$ and
$R_{c,{\bf m}}(I_{2k,2k-1})$ we have the expressions
$$
R_{c,{\bf m}}(I_{2k+1,2k})|{\bf m}_{2k}, M\rangle =
\sum _{j=1} ^k
([c+l_{j,2k}][c-l_{j,2k}-1])^{1/2}
\frac{{\tilde A}_{2k}^j ({\bf m}_{2k},M)}
{q^{l_{j,2k}}+q^{-l_{j,2k}}}
|{\bf m}^{+j}_{2k}, M_{2k}\rangle - $$
$$\qquad\qquad\qquad
-\sum _{j=1} ^k  ([c+l_{j,2k}-1][c-l_{j,2k}])^{1/2}
\frac{{\tilde A}_{2k}^j ({\bf m}^{-j}_{2k},M)}
{q^{l_{j,2k}}+q^{-l_{j,2k}}}
| m^{-j}_{2k},M\rangle , \eqno (32) $$
$$
R_{c,{\bf m}}(I_{2k,2k-1})|{\bf m}_{2k-1}, M\rangle =
\sum _{j=1} ^{k-1}
([c+l_{j,2k-1}][c-l_{j,2k-1}])^{1/2}
\frac{ {\tilde B}_{2k-1}^j ({\bf m}_{2k-1},M)}
{[2l_{j,2k-1}-1][l_{j,2k-1}]}
| {\bf m}^{+j}_{2k-1},M\rangle - $$
$$
- \sum _{j=1} ^{k-1} ([c+l_{j,2k-1}-1][c-l_{j,2k-1}+1])^{1/2}
\frac{ {\tilde B}_{2k-1}^j ({\bf m}^{-j}_{2k-1},M)}
{[2l_{j,2k-1}-1][l_{j,2k-1}-1]}
| {\bf m}^{-j}_{2k-1},M\rangle + $$
$$
+{\rm i}[c]{\tilde C}_{2k-1}({\bf m}_{2k-1},M) |{\bf m}_{2k-1},M\rangle ,
\eqno (33)
$$
where ${\tilde A}_{2k}^j ({\bf m}_{2k},M)$,
${\tilde B}_{2k-1}^j ({\bf m}_{2k-1},M)$ and
${\tilde C}_{2k-1}({\bf m}_{2k-1},M)$ are such as in (27) and (28).
These formulas were given in [20] (see also [8]).
\medskip

\noindent
{\bf Theorem 2.} {\it The representation $R_{c,{\bf m}}$ of
$U'_q({\rm so}_{2k,1})$ is irreducible if and only if
$c$ is not integer (resp. half-integer) if $l_{j,2k+1}$, $j=2,3,\cdots ,k$,
are integers (resp. half-integers) or one of the numbers $c$, $1-c$ coincides
with one of the numbers $l_{j,2k+1}$, $j=2,3,\cdots ,k$.
The representation $R_{c,{\bf m}}$ of
$U'_q({\rm so}_{2k-1,1})$ is irreducible if and only if
$c$ is not integer (resp. half-integer) if $l_{j,2k}$, $j=2,3,\cdots ,k$,
are integers (resp. half-integers) or $|c|$
coincides with one of the numbers $l_{j,2k}$, $j=2,3,\cdots ,k$,
or $|c|<|l_{k,2k}|$.}
\medskip

The reducible representations $R_{c,{\bf m}}$ can be analysed as in the
case of the principal nonunitary series of the algebra ${\rm so}_{n,1}$
(see, for example, [19]). It will be made in a separate paper.
\medskip

\centerline{\bf 9. Nonclassical type representations of
$U'_q({\rm so}_{n,1})$}
\medskip

Now we describe irreducible representations of nonclassical type (that is,
representations $R$ for which there exists no limit $q\to 1$
for the operators $R(I_{i,i-1})$). These representations are given by
$\epsilon :=(\epsilon _2,\epsilon _3,\cdots ,\epsilon _{n+1})$, by
a complex
parameter $c$ and by numbers ${\bf m}=(m_{2,n+1},m_{3,n+1},\cdots$,
$m_{\{ (n+1)/2\} ,n+1})$, $m_{2,n+1}\ge m_{3,n+2} \ge \cdots \ge
m_{\{ (n+1)/2\} ,n+1}\ge 1/2$,
describing irreducible representations of
the nonclassical type of the subalgebra $U'_q({\rm so}_{n-1})$ (see
section 4). We denote the corresponding representations of
$U_q({\rm so}_{n,1})$ by $R_{\epsilon ,c ,{\bf m}}$.

In order to describe the space of the representation
$R_{\epsilon,c,{\bf m}}$ we note that
$$
R_{\epsilon,\lambda ,{\bf m}}\downarrow U'_q({\rm so}_{n}) =
\bigoplus _{{\bf m}_n}
T_{\epsilon ',{\bf m}_n},\ \ \ {\bf m}_n=(m_{1,n},\cdots ,m_{\{ n/2 \} ,n}) ,
\eqno (34)
$$
where $\epsilon '=(\epsilon _2,\cdots ,\epsilon _{n})$,
the summation is over all irreducible nonclassical type
representations $T_{\epsilon ',{\bf m}_n}$
of the subalgebra $U'_q({\rm so}_{n})$ for which
the components of ${\bf m}_{n}$ satisfy the
"betweenness" conditions
$$
m_{1,2k+1}\ge m_{1,2k}\ge m_{2,2k+1} \ge m_{2,2k} \ge ...
\ge m_{k,2k+1} \ge m_{k,2k} \ge 1/2 \ \ {\rm if}\ \ n=2k
                                                     \eqno(35)
$$
$$
m_{1,2k}\ge m_{1,2k-1}\ge m_{2,2k} \ge m_{2,2k-1} \ge ...
\ge m_{k-1,2k-1} \ge m_{k,2k} \ \ {\rm if}\ \ n=2k-1.
                                                      \eqno(36)
$$
The carrier space ${\hat {\cal H}}_{\epsilon,\bf m}$ of the representation
$R_{\epsilon,c,{\bf m}}$ decomposes as
${\hat {\cal H}}_{\epsilon, \bf m}=\bigoplus _{{\bf m}_n}
{\cal H}_{\epsilon,{\bf m}_n}$, where
the summation is such as in (34) and ${\cal H}_{\epsilon',{\bf m}_n}$ are the
subspaces, where the representations $T_{\epsilon',{\bf m}_n}$ of
$U'_q({\rm so}_{n})$ are realized. We choose the basis in every subspace
${\cal H}_{\epsilon,{\bf m}_n}$ as in section 4. The set of all these
bases gives a basis of the space ${\hat {\cal H}}_{\epsilon,\bf m}$.
We denote the basis elements by $| {\bf m}_n,M\rangle$, where $M$ are
the corresponding tableaux. The numbers $m_{ij}$ from
$| {\bf m}_n,M\rangle$ determine the numbers $l_{ij}$ as in section 3.
The numbers $m_{i,n+1}$ determine the numbers
$l_{i,n+1}$ as in section 8.
The operators $R_{\epsilon,c, {\bf m}}(I_{i,i-1})$, $i\le n$, are given by
formulas of the nonclassical type representations of the algebra
$U'_q({\rm so}_{n})$ given in section 4.
For the operators $R_{\epsilon,c, {\bf m}}(I_{2k+1,2k})$ if $n=2k$ and
$R_{\epsilon,c, {\bf m}}(T_{2k,2k-1})$ if $n=2k-1$ we have the
expressions
$$\leqno
R_{\epsilon,c, {\bf m}}(I_{2k,2k-1})|{\bf m}_{2k-1}, M\rangle =
$$
$$
=\sum _{j=1} ^{k-1}
([c+l_{j,2k-1}][c-l_{j,2k}])^{1/2}
\frac{{\tilde B}_{2k-1}^j ({\bf m}_{2k-1},M)}{[2l_{j,2k-1}-1][l_{j,2k-1}]_+}
| {\bf m}^{+j}_{2k-1},M\rangle - $$
$$
- \sum _{j=1} ^{k-1}
([c+l_{j,2k-1}+1][c-l_{j,2k}+1])^{1/2}
\frac{ {\tilde B}_{2k-1}^j ({\bf m}^{-j}_{2k-1},M)}
{[2l_{j,2k-1}-1][l_{j,2k-1}-1]_+}
| {\bf m}^{-j}_{2k-1},M\rangle +    $$
$$
+ \epsilon _{2k}
[c]_+ {\hat C}_{2k-1}({\bf m}_{2k-1},M) |{\bf m}_{2k-1},M\rangle ,
$$
$$
R_{\epsilon,c,{\bf m}}(I_{2k+1,2k})|{\bf m}_{2k}, M\rangle =
\delta_{m_{k,2k},1/2}\,[c-1/2]\, \frac{\epsilon _{2k+1}}
{q^{1/2}-q^{-1/2}} D_{2k}(m_{2k},M)
|{\bf m}_{2k}, M\rangle +     $$
$$
+ \sum _{j=1} ^k
([c+l_{j,2k}][c-l_{j,2k}-1])^{1/2}
\frac {{\tilde A}_{2k}^j ({\bf m}_{2k},M)}{q^{l_{j,2k}}-q^{-l_{j,2k}}}
|{\bf m}^{+j}_{2k}, M\rangle -       $$
$$
-
\sum _{j=1} ^k ([c+l_{j,2k}-1][c-l_{j,2k}])^{1/2} \frac
{{\tilde A}_{2k}^j ({\bf m}^{-j}_{2k},M)}{q^{l_{j,2k}}-q^{-l_{j,2k}}}
| {\bf m}^{-j}_{2k},M\rangle , $$
where the summation in the last sum must be from 1 to $k-1$ if
$m_{k,2k}=1/2$, and
${\tilde A}_{2k}^j $, ${\tilde B}_{2k-1}^j $
are such as in (27) and (28), and
$$
{\hat C}_{2k-1}(m_{2k-1},M)=
\frac {\prod _{s=2}^k [l_{s,2k}]_+\prod _{s=1}^{k-1} [l_{s,2k-2}]_+}
{\prod _{s=1}^{k-1} [l_{s,2k-1}]_+[l_{s,2k-1}-1]_+} ,
$$
$$
D_{2k}(m_{2k},M) =
\frac{\prod_{i=2}^k
[l_{i,2k+1}-\frac 12 ] \prod_{i=1}^{k-1} [l_{i,2k-1}-\frac 12 ] }
{\prod_{i=1}^{k-1}
[l_{i,2k}+\frac 12 ] [l_{i,2k}-\frac 12 ] } .
$$
{\bf Theorem 3.} {\it The representation $R_{\epsilon ,c,{\bf m}}$ of
$U'_q({\rm so}_{2k,1})$ is irreducible if and only if
$c$ is not half-integer or one of the numbers $c$, $1-c$ coincides
with one of the numbers $l_{j,2k+1}$, $j=2,3,\cdots ,k$.
The representation $R_{\epsilon ,c,{\bf m}}$ of
$U'_q({\rm so}_{2k-1,1})$ is irreducible if and only if
$c$ is not half-integer or $|c|$
coincides with one of the numbers $l_{j,2k}$, $j=2,3,\cdots ,k$,
or $|c|<l_{k,2k}$.}
\medskip

This theorem will be proved in a separate paper. There will be also studied
reducible representations $R_{\epsilon ,c,{\bf m}}$.

\end{document}